\documentclass[10pt]{article}

\usepackage{graphicx}                           
\usepackage{color}                              
\usepackage{indentfirst}                        
\usepackage{amsmath,amssymb,amsfonts,amsthm,bm} 
\usepackage{cases}                              
\usepackage{cite}
\usepackage{mathrsfs}

\newtheorem{theorem}{Theorem}[section]
\newtheorem{definition}{Definition}[section]
\newtheorem{lemma}{Lemma}[section]

\newtheorem{example}{Example}[section]
\newtheorem{remark}{Remark}[section]

\theoremstyle{remark}
\theoremstyle{remark}
\begin{document}
\title{Gerghaty type results via simulation and $\mathcal{C}$-class functions with application}
\author{Azhar Hussain\thanks{
Department of Mathematics, University of Sargodha, Sargodha-40100, Pakistan. Email: hafiziqbal30@yahoo.com}, Muhammad Ishfaq\thanks{
Department of Mathematics, University of Sargodha, Sargodha-40100, Pakistan. Email: ishfaqahmad5632@gmail.com,} Tanzeela Kanwal \thanks{Department of Mathematics, University of Sargodha, Sargodha-40100, Pakistan. Email: tanzeelakanwal16@gmail.com,} Stojan Radenovi\'{c}\thanks{University of Belgrade, Faculty of Mechanical Engineering, Kraljice Marije
16, 11120 Beograd 35, Serbia, Department of Mathematics, College of Science, King Saud University,
Saudi Arabia, Email: radens@beotel.rs} }
\date{}
\maketitle

\noindent\textbf{ Abstract}: In this paper we study the notion of Gerghaty type contractive mapping via simulation function along with $\mathcal{C}$-class functions and prove the existence of several fixed point results in ordinary and partially ordered metric spaces. An example is given to show the validity of our results given herein. Moreover, existence of solution of two-point boundary value second order nonlinear differential equation is obtain.  \\



\section{Introduction}
The Banach contraction principle \cite{banach} is one of the fundamental result
in metric fixed point theory. Because of its importance in nonlinear analysis, number of authors
have improved, generalized and extended this basic result either by defining a new contractive
mapping in the context of a complete metric space or by investigating the existing
contractive mappings in various abstract spaces (see, e.g., \cite{agar, Ciric, cho, ZM, niet, amin} and references therein).

In particularly, Geraghty \cite{ger} consider an auxiliary function and generalized the Banach contraction in the frame work of complete metric space. Later on, Amini-Harandi and Emami \cite{amin} obtained similar results in the setting of partially ordered metric spaces. Using the concept of Samet \cite{samet}, Cho et al. \cite{cho} generalized Geraghty contraction to $\alpha$-Geraghty contraction and prove a fixed point theorem for such contraction.
On the other hand, Khojasteh {\it et al.} \cite{stojan2} introduced the notion of $\mathcal{Z}$-contraction by using a function called simulation function and proved a version of Banach contraction principle.

\section{Preliminaries}
In this section we present some basic notions and results from the literature: \\
We denote by $\mathcal{F}$ the class of all functions $ \beta : [0,\infty) \rightarrow [0,1)$ satisfying
$ \beta (t_{n}) \to 1$, implies $t_{n} \rightarrow 0$ as $n \to \infty$.
\begin{definition}\cite{ger}
Let $(X,d)$ be a metric space. A map $T:X\rightarrow X$ is called Geraghty contraction if there exists $\beta\in \mathcal{F}$ such that for all $x, y\in X$,
$$d(Tx,Ty)\leq\beta(d(x,y))d(x,y).$$
\end{definition}
\begin{theorem}\cite{ger}
Let $(X,d)$ be a complete metric space. Mapping $T:X \rightarrow X$ is Geraghty contraction.
Then $T$ has a fixed point $x\in X$, and $\{T^{n}x_{1}\}$ converges to $x$.
\end{theorem}
In 2015, Khojasteh {\it et al.} \cite{stojan2} introduced simulation function $\zeta : [0, \infty) \times  [0, \infty) \rightarrow \mathbb{R}$, satisfying the following assertions:
\begin{description}
  \item[$(\zeta_{1})$] $\zeta(0, 0)=0$;
  \item[$(\zeta_{2})$]  $\zeta(t, s)<s-t$ for all $t, s>0$;
  \item[$(\zeta_{3})$]  If $\{t_{n}\}, \{s_{n}\}$ are sequences in $(0, \infty)$ such that $\lim\limits_{n\rightarrow \infty}t_{n}=\lim\limits_{n\to \infty}s_{n}>0$ then $$\lim\limits_{n\rightarrow \infty}sup~\zeta(t_{n}, s_{n})<0$$
\end{description}
\begin{definition}\cite{stojan2}
Let $(X, d)$ be a metric space, $T: X \rightarrow X$ a mapping and $\zeta$ a simulation function. Then $T$ is called a $\mathcal{Z}$-contraction with respect to $\zeta$ if it satisfies
\begin{equation}\label{Eq01}
\zeta(d(Tu, Tv), d(u, v)) \geq 0 ~~\text{for all}~ u, v \in X.
\end{equation}
\end{definition}
\begin{theorem}\label{T1}\cite{stojan2}
Let $(X, d)$ be a complete metric space and $T : X \rightarrow X$ be a $\mathcal{Z}$-contraction with respect to $\zeta$. Then $T$ has a unique fixed point $u\in X$ and for every $x_{0} \in X$, the Picard sequence $\{x_{n}\}$ where $x_{n}=Tx_{n-1}$ for all $n \in \mathbb{N}$ converges to this fixed point of $T$.
\end{theorem}
\begin{example}\cite{stojan2}
Let $ \zeta_{i} : [0, \infty) \times [0, \infty)\rightarrow \mathbb{R},~~i=1, 2, 3$ be defined by
\begin{description}
  \item[$(i)$] $\zeta_{1}(t, s)=\lambda s-t$, where $\lambda \in (0, 1)$;
  \item[$(ii)$]  $\zeta_{2}(t, s)= s\varphi(s)-t$, where $\varphi:[0, \infty)\rightarrow [0, 1)$ is a mapping such that $\lim\limits_{t\rightarrow r^{+}}sup~\psi(t)<1$ for all $r>0$;
  \item[$(iii)$]  $\zeta_{3}=s-\psi(s)-t$, where $\psi:[0, \infty)\rightarrow [0, \infty)$ is a continuous function such that $\psi(t)=0$ if and only if $t=0.$
\end{description}
Then $\zeta_{i}$ for $i=1, 2, 3$ are simulation functions.
\end{example}
Rold$\acute{a}$n-L$\acute{o}$pez-de-Hierro {\it et al.} \cite{Roldan} modified the notion of a simulation function by replacing $(\zeta_{3})$ by $(\zeta'_{3})$,
\begin{description}
\item[($\zeta'_{3}$)]: if $\{t_{n}\}, \{s_{n}\}$ are sequences in $(0, \infty)$ such that $\lim_{n\rightarrow \infty}t_{n}=\lim\limits_{n\rightarrow \infty}s_{n}>0$ and $t_{n}<s_{n}$, then $$\lim\limits_{n\rightarrow \infty}sup\zeta(t_{n}, s_{n})<0.$$
\end{description}
The function $\zeta: [0, \infty) \times [0, \infty) \rightarrow \mathbb{R}$ satisfying $(\zeta_{1}-\zeta{2})$ and $(\zeta'_{3})$ is called simulation function in the sense of Rold$\acute{a}$n-L$\acute{o}$pez-de-Hierro.

\begin{definition}\cite{arslan}\label{D4.2}
A mapping $\mathcal{G}:[0,+\infty)^{2}\rightarrow \mathbb{R}$ is called a
$\mathcal{C}$-class function if it is continuous and satisfies the following
conditions:
\begin{description}
  \item [(1)]  $\mathcal{G}(s, t)\leq s$;
  \item [(2)]  $\mathcal{G}(s, t)= s$ implies that either $s=0$ or $t=0$,
      for all $s, t \in [0, +\infty)$.
\end{description}
\end{definition}
\begin{definition}\cite{arslan2}\label{D4.3}
A mapping $\mathcal{G}:[0,+\infty)^{2}\rightarrow\mathbb{R}$ has the
property $\mathcal{C}_{\mathcal{G}}$, if there exists and $\mathcal{C}_{\mathcal{G}}\geq 0$ such
that
\begin{description}
  \item [(1)]  $\mathcal{G}(s, t)> \mathcal{C}_{\mathcal{G}}$ implies $s>t$;
  \item [(2)]  $\mathcal{G}(s, t)\leq \mathcal{C}_{\mathcal{G}}$, for all $t\in[0,
      +\infty)$.
\end{description}
\end{definition}
Some examples of $\mathcal{C}$-class functions that have property $\mathcal{C}_{\mathcal{G}}$ are
as follows:
\begin{description}
  \item [(a)] $\mathcal{G}(s, t)= s-t,~ \mathcal{C}_{\mathcal{G}}=r, r\in[0,
      +\infty)$;
  \item [(b)] $\mathcal{G}(s, t)= s-\frac{(2+t)t}{(1+t)},
      \mathcal{C}_{\mathcal{G}}=0$;
  \item[(c)] $\mathcal{G}(s, t)= \frac{s}{1+kt}, k\geq1,
      \mathcal{C}_{\mathcal{G}}=\frac{r}{1+k}, r\geq 2$.
\end{description}
For more examples of $\mathcal{C}$-class functions that have property $\mathcal{C}_{\mathcal{G}}$
see \cite{arslan3, ZM, arslan2}.
\begin{definition}\cite{arslan2}\label{D4.5}
A $\mathcal{C}_{\mathcal{G}}$ simulation function is a mapping
$\mathcal{G}:[0,+\infty)^{2}\rightarrow\mathcal{R}$ satisfying the following
conditions:
\begin{description}
  \item [(1)]  $\zeta(t, s)< \mathcal{G}(s, t)$ for all $t, s>0$, where
      $\mathcal{G}:[0, +\infty)^{2}\rightarrow \mathbb{R}$ is a $\mathcal{C}$-class
      function;
  \item [(2)] if $\{t_{n}\}, \{s_{n}\}$ are sequences in $(0, +\infty)$
      such that $\lim\limits_{n\rightarrow \infty}t_{n}=
      \lim\limits_{n\rightarrow \infty}s_{n}>0 $, and $t_{n}<s_{n}$, then
      $\lim\limits_{n\rightarrow\infty}\sup\zeta(t_{n},
      s_{n})<\mathcal{C}_{\mathcal{G}}$.
\end{description}
\end{definition}
Some examples of simulation functions and $\mathcal{C}_{\mathcal{G}}$-simulation
functions are:
\begin{description}
  \item [(1)]  $\zeta(t, s)= \frac{s}{s+1}-t$ for all $t, s>0$.
  \item [(2)] $\zeta(t, s)= s-\phi(s)-t$ for all $t, s>0$, where
      $\phi:[0, +\infty)\rightarrow [0, +\infty)$ is a lower semi
      continuous function and $\phi(t)=0$ if and only if $t=0$.
\end{description}
For more examples of simulation functions and $\mathcal{C}_{\mathcal{G}}$-simulation
functions see \cite{arslan3, Roldan, stojan2, arslan2, vetro, wang}.

\begin{definition}\cite{samet}
Let $T:A\rightarrow B$ be a map and $\alpha :X \times X \rightarrow \mathbb{R}$ be a function. Then $T$ is said to be $\alpha$-admissible if $\alpha(x,y)\geq 1$
 implies $\alpha(Tx,Ty)\geq 1.$
\end{definition}
\begin{definition}\cite{karap}
An $\alpha$-admissible map $T$ is said to be triangular $\alpha$-admissible if $\alpha (x,z)\geq 1$ and $\alpha (z,y)\geq 1$  implies $\alpha(x,y)\geq 1$
\end{definition}

Cho et al. \cite{cho} generalized the concept of Geraghty contraction to $\alpha$-Geraghty contraction and prove the fixed point theorem for such contraction.

\begin{definition}\cite{cho}
Let $(X,d)$ be a metric space, and let $\alpha : X \times X \rightarrow \mathbb{R}$ be a function. A map $T:X\rightarrow X$ is called $\alpha$-Geraghty contraction if there exists $\beta\in \mathcal{F}$ such that for all $x, y\in X$,
$$\alpha (x,y)d(Tx,Ty)\leq\beta(d(x,y))d(x,y).$$
\end{definition}
\begin{theorem}\cite{cho}
Let $(X,d)$ be a complete metric space, $\alpha: X\times X \rightarrow \mathbb{R}$ be a function. Define a map $T:X \rightarrow X$ satisfying the following conditions:
\begin{enumerate}
  \item $T$ is continuous $\alpha$-Geraghty contraction;
  \item $T$ be a triangular $\alpha$-admissible;
  \item there exists $x_{1}\in X$ such that $\alpha(x_{1}, Tx_{1})\geq 1$;
\end{enumerate}
Then $T$ has a fixed point $x\in X$, and $\{T^{n}x_{1}\}$ converges to $x$.
\end{theorem}
\begin{lemma}\cite{karap}\label{L3.1}
Let $T:X\rightarrow X$ be a triangular $\alpha$-admissible map. Assume that
there exists $x_{1} \subset X$ such that $ \alpha(x_{1}, Tx_{1})\geq 1
$.Define a sequence $\{x_{n}\}$ by $x_{n+1}= Tx_{n}$. Then we have
$\alpha(x_{n}, x_{m})\geq 1$ for all $m, n \in \mathbb{N}$ with $n < m$.
\end{lemma}

\begin{lemma}\cite{stojan}\label{L4.2}
Let $(X, d)$ be a metric space and let $\{x_{n}\}$ be a sequence in $X$ such
that
\begin{equation}
\lim\limits_{n\rightarrow\infty}d(x_{n}, x_{n+1})=0.
\end{equation}
If $\{x_{n}\}$ is not a Cauchy sequence in $X$, then there exists
$\varepsilon >0$ and two sequences $x_{m(k)}$ and $x_{n(k)}$ of positive
integers such that $x_{n(k)}>x_{m(k)}>k$ and the following sequences tend to
$\varepsilon^{+}$ when $k\rightarrow\infty$:
$$d(x_{m(k)},x_{n(k)}), d(x_{m(k)},x_{n(k)+1}),d(x_{m(k)-1},x_{n(k)}),$$ $$d(x_{m(k)-1},x_{n(k)+1}),d(x_{m(k)+1},x_{n(k)+1}).$$
\end{lemma}
Motivated by the above results, we introduce the notion of Gerghaty type $\mathcal{Z}_{(\alpha, \mathcal{G})}$-contraction and prove some fixed point results in metric and partially ordered metric spaces. An example to prove the validity and application to nonlinear differential equation for the usability of our results is presented.

\section{Fixed point results in usual metric space}
We begin with the following notion:

\begin{definition}\label{D3.3}
Let $(X, d)$ be a metric space and $\alpha: X\times
X\rightarrow [0, \infty)$ be a function. A mapping $T:X\rightarrow X$ is called a
$\mathcal{Z}_{(\alpha, \mathcal{G})}$-Geraghty contraction if there exists $\beta
\in \mathcal{F}$ such that for all $x, y \in X$
\begin{equation}\label{E3.1}
\zeta(\alpha(x, y)d(Tx, Ty), \beta(M(x, y)) M(x, y))\geq \mathcal{C}_{\mathcal{G}}
\end{equation}
where
$$ M(x, y)= \max\{d(x, y), d(x, Tx), d(y, Ty)\}.$$
\end{definition}
\begin{remark}\label{R3.1}
Since the functions belonging to $\mathcal{F}$ are strictly smaller than 1,
$(\ref{E3.1})$ implies that
$$ d(Tx, Ty) < M(x, y)$$ for any $x, y \in X$ with $x \neq y$ and for $\alpha(x, y)\geq 1$, $\zeta(t,s)<\mathcal{G}(s, t)=s-t$.
\end{remark}
\begin{theorem}\label{Th3.1}
Let $(X, d)$ be a complete metric space, $\alpha:X\times X\rightarrow [0,
\infty)$ and $T: X \rightarrow X$ be two functions. Suppose that the
following conditions are satisfied:
\begin{description}
  \item [(1)]  $T$ is continuous $Z_{(\alpha, \mathcal{G})}$-Geraghty contraction;
  \item [(2)] $T$ is triangular $\alpha$-admissible;
  \item [(3)] there exists $x_{1} \in X $ such that $ \alpha (x_{1},
      Tx_{1})\geq 1$;
  \item [(4)] $T$ is continuous.
\end{description}
Then $T$ has a fixed point $x^{*} \in X $ and $T$ is a Picard operator that
is, ${T^{n}x_{1}}$ converges to $x ^{*}$.
\begin{proof}
Let $x_{1}\in X$ be such that $\alpha (x_{1}, Tx_{1})\geq 1$. Define a
sequence $\{x_{n}\}\subset X $ by $x_{n+1}= Tx_{n} $ for $n \in \mathbb{N}$.
If $x _{n_{0}}= x_{n_{0}+ 1}$ for some $n_{0}\in\mathbb{N}$, then $x_{n_{0}}$
is a fixed point of $T$ and hence the proof is completed. Thus, we assume that $x_{n}\neq x_{n+1}$ for all $n\in \mathbb{N}$. By
Lemma \ref{L3.1}, we have
\begin{equation}\label{E3.2}
\alpha(x_{n}, x_{n+1}) \geq 1
\end{equation}
for all $n\in \mathbb{N}$. Then
\begin{eqnarray}\label{Eq3.1}
d(x_{n+1}, x_{n+2}) &=&d(Tx_{n}, Tx_{n+1})\nonumber \\
 &\leq&\alpha(x_{n}, x_{n+1}) d(Tx_{n}, Tx_{n+1}).
\end{eqnarray}
Since $T$ is a $\mathcal{Z}_{(\alpha, \mathcal{G})}$-Geraghty contraction, we have
\begin{eqnarray*}
\mathcal{C}_{\mathcal{G}}&\leq& \zeta (\alpha(x_{n}, x_{n+1}) d(Tx_{n}, Tx_{n+1}), \beta (M(x_{n}, x_{n+1}))M(x_{n}, x_{n+1}))\\
&<&\mathcal{G}(\beta (M(x_{n}, x_{n+1}))M(x_{n}, x_{n+1}), \alpha(x_{n}, x_{n+1}) d(Tx_{n}, Tx_{n+1})).
\end{eqnarray*}
Using $(\mathcal{G}_{1})$, we obtain
\begin{equation}\label{E3.3}
\alpha(x_{n}, x_{n+1}) d(Tx_{n}, Tx_{n+1})< \beta (M(x_{n}, x_{n+1})) M(x_{n}, x_{n+1}).
\end{equation}
From (\ref{Eq3.1}) and (\ref{E3.3}), we have
\begin{equation}\label{E3.4}
d(x_{n+1}, x_{n+2})< \beta (M(x_{n}, x_{n+1})) M(x_{n}, x_{n+1})
\end{equation}
for all $n \in \mathbb{N}$, where
\begin{eqnarray*}
 M(x_{n}, x_{n+1})&=&\max\{ d(x_{n}, x_{n+1}), d(x_{n},Tx_{n}), d(x_{n+1}, Tx_{n+1})\}\\
&=&\max\{d(x_{n}, x_{n+1}), d(x_{n}, x_{n+1}), d(x_{n+1}, x_{n+2})\} \\
&=&\max\{d(x_{n}, x_{n+1}), d(x_{n+1}, x_{n+2})\}.
\end{eqnarray*}
If $M(x_{n},x_{n+1})=d(x_{n+1},x_{n+2})$, then by definition of $\beta$, we have
\begin{eqnarray*}
d(x_{n+1}, x_{n+2}) &<& \beta (d(x_{n}, x_{n+1})d(x_{n}, x_{n+1}))\\
&<& d(x_{n+1}, x_{n+2})
\end{eqnarray*}
a contradiction. Thus we conclude that $M(x_{n}, x_{n+1}) = d(x_{n}, x_{n+1})$  for all $n \in \mathbb{N}$ and so the sequence $\{d(x_{n}, x_{n+1})\}$
of real numbers is decreasing and bounded below by zero. Hence there exists $r \geq 0 $ such that
$$\lim\limits_{n\rightarrow \infty} d(x_{n}, x_{n+1})= r.$$
We claim that $r=0$. Suppose on contrary that $r > 0$, then by (\ref{E3.4}) we have
\begin{eqnarray*}
\frac {d(x_{n+1}, x_{n+2})}{d(x_{n}, x_{n+1})} \leq  \beta (d(x_{n}, x_{n+1}) < 1
\end{eqnarray*}
 which yields that $\lim\limits_{ n\rightarrow \infty }\beta (d(x_{n}, x_{n+1})= 1$.
 Since $\beta \in \mathcal{F} $, we get that
 \begin{equation}\label{E3.5}
 \lim\limits_{ n\rightarrow \infty}  d(x_{n}, x_{n+1})= 0.
 \end{equation}
 We now show that $\{x_{n}\}$ is a Cauchy sequence. Suppose on contrary that
 it is not. Thus there exists $\epsilon > 0 $ such that
 for all $k > 0$, $m(k) > n(k)> k$ with the (smallest number
 satisfying the condition below)
$ d(x_{m(k)}, x_{n(k)}) \geq \epsilon $  and $ d(x_{m(k)-1}, x_{n(k)}) <
\epsilon $.
 Then we have
 \begin{eqnarray*}
\epsilon &\leq& d(x_{m(k)}, x_{n(k)})\\
&\leq&  d(x_{m(k)}, x_{m(k)-1}) + d(x_{m(k)-1}, x_{n(k)})\\
&\leq& d(x_{m(k)}, x_{m(k)-1}) +\epsilon.
 \end{eqnarray*}
Letting $k \rightarrow \infty $ in the above inequality, we have
\begin{equation}\label{E3.6}
\lim\limits_{k\to \infty}d(x_{m(k)}, x_{n(k)}) =   \epsilon.
 \end{equation}
 By using (\ref{E3.5}) and (\ref{E3.6}), we obtain
 $$ \lim\limits_{k\to \infty} d(x_{m(k)-1}, x_{n(k)-1}) =   \epsilon. $$
 By Lemma \ref{L3.1},  $\alpha(x_{m(k)-1}, x_{n(k)-1})\geq 1$, thus
\begin{eqnarray}\label{E3.6'}
d(x_{m(k)}, x_{n(k)}) &=& d(Tx_{m(k)-1}, Tx_{n(k)-1})\nonumber \\
&\leq & \alpha(x_{m(k)-1}, x_{n(k)-1})d(Tx_{m(k)-1}, Tx_{n(k)-1}).
\end{eqnarray}
Since $T$ is a $\mathcal{Z}_{(\alpha, \mathcal{G})}$-Geraghty contraction, we have
\begin{eqnarray*}
\zeta ( \alpha(x_{m(k)-1}, x_{n(k)-1})d(Tx_{m(k)-1}, Tx_{n(k)-1}), \beta (M(x_{m(k)-1}, x_{n(k)-1})) M(x_{m(k)-1}, x_{n(k)-1})) \geq \mathcal{C}_{\mathcal{G}}.
\end{eqnarray*}
This implies
\begin{eqnarray*}
&&\mathcal{C}_{\mathcal{G}}\\&\leq& \zeta (\alpha(x_{m(k)-1}, x_{n(k)-1})d(Tx_{m(k)-1}, Tx_{n(k)-1}), \beta (M(x_{m(k)-1}, x_{n(k)-1})) M(x_{m(k)-1}, x_{n(k)-1}))\\
&<&\mathcal{G}(\beta (M(x_{m(k)-1}, x_{n(k)-1})) M(x_{m(k)-1}, x_{n(k)-1}), \alpha(x_{m(k)-1}, x_{n(k)-1})d(Tx_{m(k)-1}, Tx_{n(k)-1})).
\end{eqnarray*}
Using $(\mathcal{G}_{1})$, we obtain
\begin{equation}\label{E3.6''}
\alpha(x_{m(k)-1}, x_{n(k)-1})d(Tx_{m(k)-1}, Tx_{n(k)-1}))< \beta(M(x_{m(k)-1}, x_{n(k)-1})) M(x_{m(k)-1}, x_{n(k)-1}),
\end{equation}
where
\begin{eqnarray*}
&& M((x_{m(k)-1}, x_{n(k)-1}))\\&=& \max\{ d(x_{m(k)-1}, x_{n(k)-1}), d((x_{m(k)-1},Tx_{m(k)-1}), d(x_{n(k)-1}, Tx_{n(k)-1})\}\\
&=& \max\{ d(x_{m(k)-1}, x_{n(k)-1}), d(x_{m(k)-1}, x_{m(k)}),
d(x_{n(k)-1}, x_{n(k)})\}.
\end{eqnarray*}
If $M(x_{m(k)-1}, x_{n(k)-1})=d(x_{m(k)-1}, x_{m(k)})$, we have
\begin{eqnarray*}
\alpha(x_{m(k)-1}, x_{n(k)-1})d(x_{m(k)}, x_{n(k)}) &<& \beta (d(x_{m(k)-1}, x_{m(k)})) d(x_{m(k)-1}, x_{m(k)})\\
&<&  d(x_{m(k)-1}, x_{m(k)})
\end{eqnarray*}
a contradiction. Similarly, we have contradiction when $M(x_{m(k)-1}, x_{n(k)-1})=d(x_{n(k)-1},
x_{n(k)})$. Thus we conclude that $M(x_{m(k)-1}, x_{n(k)-1})=d(x_{m(k)-1},
x_{n(k)-1})$. So
\begin{eqnarray*}
 \frac{d(x_{m(k)}, x_{n(k)})}{d(x_{m(k)-1}, x_{n(k)-1})} \leq
\beta (d(x_{m(k)-1}, x_{n(k)-1}) < 1.
\end{eqnarray*}
Letting $k\rightarrow \infty$ in above inequality, we derive that
$$\lim\limits _{k\rightarrow \infty } \beta(d(x_{m(k)-1}, x_{n(k)-1})) = 1.$$
This implies
 $$\lim\limits_{k\rightarrow \infty} d(x_{m(k)-1}, x_{n(k)-1}) = 0.$$
 Hence $\epsilon = o$, which is a contradiction. Thus we conclude that
 $\{x_{n}\}$ is a Cauchy sequence.
 It follows from completeness of $X$ that there exists $x^{*}\in X$ such that
 $$\lim\limits_{n\rightarrow \infty} x_{n}= x^{*}.$$
 Since $T$ is continuous, we get $ \lim\limits_{n\rightarrow \infty} x_{n+1} =
 Tx_{n}$
 and so $x^{*}= Tx^{*}$.
This completes the proof.
\end{proof}
\end{theorem}
\begin{theorem}\label{Th3.2}
Let $(X, d)$ be a complete metric space, $\alpha:X\times X\rightarrow[0,
\infty)$ and $T:X  \rightarrow X$ be two mappings. Suppose that the following
conditions are satisfied:
\begin{description}
  \item [(1)]  $T$ is $\mathcal{Z}_{(\alpha, \mathcal{G})}$-Geraghty contraction;
  \item [(2)] $T$ is triangular $\alpha$-admissible;
  \item [(3)] there exists $x_{1} \in X$ such that $ \alpha(x_{1},
      Tx_{1}) \geq 1$;
  \item [(4)] if $x_{n}$ is a sequence in $X$ such that $\alpha(x_{n},
      x_{n+1}) \geq 1$ for all $n$ and $x_{n} \rightarrow x^{\ast} \in X$ as $n
      \rightarrow \infty $, then there exist a subsequence $\{x_{n(k)}\}$
      of $\{x_{n}\}$ such that $\alpha(x_{n(k)}, x^{\ast}) \geq 1$ for all $k$.
\end{description}
Then $T$ has a fixed point $x^{*} \in X $ and $T$ is a Picard operator, that
is, $ {T^{n}x_{1}} $ converges to $ x ^{*}$ .
\begin{proof}
Following the arguments those given in Theorem \ref{Th3.1}, we conclude that the sequence
${x_{n}}$ defined  by $x_{n+1}= Tx_{n}$ for all  $n \geq 0 $, converges to
$x^{*} \in X$. By condition (4) we deduce that
there exists a subsequence  $\{x_{n(k)}\}$  of $\{x_{n}\}$  such that
$\alpha(x_{n(k)}, x^{*}) \geq 1$ for all $k$. Also
\begin{eqnarray}\label{Eq3.2}
d(x_{n(k)+1}, Tx^{*}) &=&d(Tx_{n(k)}, Tx^{*})\nonumber \\
&\leq&\alpha(x_{n(k)}, x^{*}) d(Tx_{n(k)}, Tx^{*}).
\end{eqnarray}
Since $T$ is $\mathcal{Z}_{(\alpha, \mathcal{G})}$-Geraghty contraction, we have
\begin{eqnarray*}
\mathcal{C}_{\mathcal{G}} \leq \zeta ( \alpha(x_{n(k)}, x^{*}) d(Tx_{n(k)}, Tx_{*}), \beta (M(x_{n(k)}, x^{*})) M(x_{n(k)}, x^{*}))\\
<\mathcal{G}(\beta (M(x_{n(k)}, x^{*})) M(x_{n(k)}, x^{*})), \alpha(x_{n(k)}, x^{*}) d(Tx_{n(k)}, Tx^{*})).
\end{eqnarray*}
By definition of $\mathcal{G}$, we get that
\begin{equation}\label{Eq3.7}
\alpha(x_{n(k)}, x^{*}) d(Tx_{n(k)}, Tx^{*})<\beta (M(x_{n(k)}, x^{*})) M(x_{n(k)}, x^{*})).
\end{equation}
From (\ref{Eq3.2}) and (\ref{Eq3.7}), we have
\begin{equation}\label{E3.8}
 d(x_{n(k)+1}, Tx^{*})< \beta (M(x_{n(k)}, x^{*})) M(x_{n(k)}, x^{*})),
\end{equation}
where
\begin{eqnarray*}
 M(x_{n(k)}, x^{*})&=&\max\{ d(x_{n(k)}, x^{*}), d(x_{n(k)},Tx_{n(k)}), d(x^{*}, Tx^{*})\}\\
&=&\max\{ d(x_{n(k)}, x^{*}), d(x_{n(k)}, x_{n(k)+1}),
d(x^{*}, Tx^{*})\}.
\end{eqnarray*}
Letting $k \rightarrow \infty $ in the above inequality, we get that
\begin{equation}\label{E3.8'}
\lim _{k \rightarrow \infty} M(x_{n(k)}, x^{*})= d(x^{*}, Tx^{*}).
 \end{equation}
 Suppose $d(x^{*}, Tx^{*}) > 0 $.
 By definition of $\beta$ and \eqref{E3.8}, we have
 $$ d(x_{n(k)+1}, Tx^{*})<M(x_{n(k)}, x^{*}).$$
Letting $k\rightarrow \infty $ in the above inequality and using (\ref{E3.8'}), we
obtain that
$$ d(x^{*}, Tx^{*}) < d(x^{*}, Tx^{*}),$$
a contradiction. Thus $ d(x^{*}, Tx^{*}) = 0$, that is, $ x^{*} = Tx^{*}.$
\end{proof}
\end{theorem}
For the uniqueness of fixed point, we consider the following hypothesis:

$(\mathcal{U})$ For all $x,y\in Fix(T)$, there exists $z\in X$ such that $\alpha(x, z)\geq 1$ and
$\alpha(y, z)\geq 1$. Here, $Fix(T)$ denotes the set of fixed points of $T$.
\begin{theorem}\label{Th3.3}
Adding condition $(\mathcal{\mathcal{U}})$ to the hypothesis of Theorem \ref{Th3.1} (resp. Theorem
\ref{Th3.2}), we obtain that $x^{*}$ is the unique fixed point of $T$.
\begin{proof}
From Theorem \ref{Th3.1} (resp. Theorem \ref{Th3.2}), we have a fixed point,
namely $x^{*} \in X$ of $T$. For uniqueness, suppose there is another fixed point of $T$, say, $y^{\ast}\in X$. Then, by assumption $(\mathcal{U})$,
there exists $z \in X$ such that
$$\alpha(x^{*}, z)\geq 1 ~ ~ and ~ ~ \alpha(y^{*}, z)\geq 1.$$
Since $T$ is a $\alpha$-admissible, we have
$$\alpha(x^{*}, T^{n}z)\geq 1 ~ ~ and ~ ~ \alpha(y^{*}, T^{n}z)\geq 1 $$
for all $n$. Hence we have
\begin{eqnarray}\label{Eq3.8''}
d(x^{*}, T^{n}z) &=&d(Tx^{*}, TT^{n-1}z)\\
&\leq&\alpha(x^{*}, T^{n-1}z)d(Tx^{*}, TT^{n-1}z).
\end{eqnarray}
Since $T$ is a $\mathcal{Z}_{(\alpha, \mathcal{G})}$-Geraghty contraction, we have
\begin{eqnarray*}
\mathcal{C}_{\mathcal{G}}\leq \zeta ( \alpha(x^{*}, T^{n-1}z) d(Tx^{*}, TT^{n-1}z), \beta (M(x^{*}, T^{n-1}z)) M(x^{*}, T^{n-1}z)\\
<\mathcal{G}(\beta (M(x^{*}, T^{n-1}z)) M(x^{*}, T^{n-1}z) ,\alpha(x^{*},  T^{n-1}z) d(Tx^{*)}, TT^{n-1}z)).
\end{eqnarray*}
By definition of $\mathcal{G}$, we have
\begin{eqnarray*}
\alpha(x^{*}, T^{n-1}z) d(Tx^{*}, TT^{n-1}z)< \beta (M(x^{*}, T^{n-1}z)) M(x^{*}, T^{n-1}z),
\end{eqnarray*}
where
\begin{eqnarray*}
  M(x^{*}, T^{n-1}z)&=&\max \{d(x^{*}, T^{n-1}z), d(x^{*}, Tx^{*}), d(T^{n-1}z, TT^{n-1}z)\}\\
&=& d(x^{*}, T^{n-1}z).
\end{eqnarray*}
Hence we have
\begin{eqnarray}\label{Eq3.8'''}
\alpha(x^{*}, T^{n-1}z) d(Tx^{*}, TT^{n-1}z) &<& \beta (d(x^{*}, T^{n-1}z)) d(x^{*}, T^{n-1}z).
\end{eqnarray}
Inequality \eqref{Eq3.8''} together with \eqref{Eq3.8'''} gives
\begin{eqnarray}\label{Eq3.9'}
d(x^{*}, T^{n}z) &<& \beta (d(x^{*}, T^{n-1}z)) d(x^{*}, T^{n-1}z).
\end{eqnarray}
By definition of $\beta$, \eqref{Eq3.9'} gives
$$d(x^{*}, T^{n}z)< d(x^{*}, T^{n-1}z)$$
for all $n \in \mathbb{N}$. Thus the sequence ${d(x^{*}, T^{n}z)}$ is non
increasing, and so there exists $ u\geq0$ such that
$\lim\limits_{n\rightarrow \infty}d(x^{*}, T^{n}z)=u$. From \eqref{Eq3.9'}, we have
$$\frac {d(x^{*}, T^{n}z)}{d(x^{*}, T^{n-1}z)} \leq \beta (d(x^{*}, T^{n-1}z))$$
and so $ \lim\limits_{n\rightarrow \infty} \beta (d(x^{*} ,T^{n-1}z)) =1 $.
Consequently, we have $ \lim\limits_{n\rightarrow \infty}(d(x^{*}, T^{n-1}z))
=0$, and hence $ \lim\limits_{n\rightarrow \infty} T^{n}z= x^{*}$. Similarly,
we can find that $ \lim\limits_{n\rightarrow \infty} T^{n}z= y^{*}$. By uniqueness of limit, we obtain
$x^{*}= y^{*}.$
\end{proof}
\end{theorem}
\begin{example}\label{Ex3.1}
Let $X=[0, \infty)$ and $d:X \times X\rightarrow \mathbb{R}$ be defined by
$d(x, y)= |x- y|$ for all $x, y \in X$. Let $\zeta(t, s)= \frac{8}{9}s -t,
\mathcal{G}(s, t)=s-t$ for all $s, t\in[0, \infty)$, $\mathcal{C}(\mathcal{G})=0$ and
$\beta(t)=\frac{1}{1+t}$ for all $t\geq 0$. Then it is clear that $\beta \in
\mathcal{F}$. We define $T: X\rightarrow X$ by
$$Tx= \begin{cases}
\frac {1}{3}x  & ~ if ~ 0 \leq x \leq 1, \\
3x & otherwise.
\end{cases}$$
and $\alpha: X\times X \rightarrow [0, \infty)$ by
$$\alpha(x, y) = \begin{cases}
1 & ~ if ~ 0 \leq x, y \leq 1, \\
0 & otherwise.
\end{cases}$$
Clearly, $T$ is continuous and condition $(3)$ of Theorem \ref{Th3.1} is satisfied with $x_{1}= 1$.
Let $x, y \in X$ such that $\alpha(x, y)\geq 1$. Then $x,
y \in [0, 1]$, so $Tx, Ty \in [0, 1]$ and thus $\alpha(Tx, Ty)=
1$. Hence $T$ is $\alpha$-admissible. Further, if $z=Ty$, then $\alpha(y,
z)\geq 1$, this implies $\alpha(x, z)\geq 1$. So $T$ is triangular $\alpha$-
admissible, hence condition $(2)$ of Theorem \ref{Th3.1} is satisfied. Finally, if $0 \leq x, y \leq 1 $,
then $\alpha(x, y)=1$, and we have
\begin{eqnarray*}
 \zeta(\alpha(x, y)d(Tx, Ty), \beta(M(x, y))M(x, y))&=&\frac{8}{9}\beta(M(x, y))M(x, y)- \alpha(x, y)d(Tx, Ty)\\
&=&\frac{8(M(x, y))}{9(1+M(x, y))}-d(Tx, Ty),
\end{eqnarray*}
where
\begin{eqnarray*}
M(x, y)&=&\max\{d(x, y), d(x, Tx), d(y, Ty)\},
\end{eqnarray*}
for all $x, y \in [0, 1]$. \\
{\bf Case-I:} If $M(x, y)=d(x, y)$, then
\begin{eqnarray*}
\zeta(\alpha(x, y)d(Tx, Ty), \beta(M(x, y))M(x, y)) &=& \frac{8d(x, y)}{9(1+d(x, y))}-d(\frac{x}{3}, \frac{y}{3})\\
&=&\frac{8|x-y|}{9(1+|x-y|)}-|\frac{x}{3}-\frac{y}{3}|\\
&=&\frac{8|x-y|}{9(1+|x-y|)}-\frac{|x-y|}{3}\\
&=&\frac{8|x-y|-3(1+|x-y|)|x-y|}{9(1+|x-y|)}\\
&=&\frac{|x-y|(5-3|x-y|)}{9(1+|x-y|)}\\
&\geq & 0.
\end{eqnarray*}
{\bf Case-II:} If $M(x, y)=d(x, Tx)$, then
\begin{eqnarray*}
\zeta(\alpha(x, y)d(Tx, Ty), \beta(M(x, y))M(x, y)) &=& \frac{8d(x, Tx)}{9(1+d(x, Tx))}-d(\frac{x}{3}, \frac{y}{3})\\
&=&\frac{8|x-\frac{x}{3}|}{9(1+|x-\frac{x}{3}|)}-|\frac{x}{3}-\frac{y}{3}|\\
&=&\frac{16|x|}{9(3+2|x|)}-\frac{|x-y|}{3}\\
&\geq& \frac{16|x|}{9(3+2|x|)}-\frac{2|x|}{9}\\
&=&\frac{|x|(10-4|x|)}{9(3+2|x|)}\\
&\geq& 0.
\end{eqnarray*}
Similarly, if $M(x, y)=d(y, Ty)$, we have
$$\zeta(\alpha(x, y)d(Tx, Ty), \beta(d(y, Ty))d(y, Ty))\geq 0.$$
Hence for $0\leq x, y \leq 1$, $T$ is a generalized
$\mathcal{Z}_{(\alpha, \mathcal{G})}$-Geraghty contraction. In either case $\alpha(x, y)= 0$ and $T$ is a $\mathcal{Z}_{(\alpha,
\mathcal{G})}$-Geraghty contraction. Thus all the hypothesis of Theorem
\ref{Th3.1} are satisfied and $T$ has a fixed point $x^{*} = 0$.
\end{example}
\section{Fixed point results in partially ordered metric space}
Let $(X, d, \preceq)$ be a partially ordered metric space. Many authors has proved the existence of fixed point results in the frame work of partially order metric spaces (see for example \cite{ran, agar, niet}). In this section, we obtain some new fixed point results in partially order metric spaces, as an application of our results given in above section.
\begin{definition}\label{D3.4}
Let $(X, d, \preceq)$ be a complete partially ordered metric space and let $x
\preceq y$ for all $x, y \in X$. A map $T:X\rightarrow X$ is called
$\mathcal{Z}_{\mathcal{G}}$-Geraghty contraction if there exists $\beta\in\mathcal{F}$
such that for all $x, y \in X$
$$\zeta(d(Tx, Ty), \beta(M(x, y))M(x, y))\geq \mathcal{C}_{\mathcal{G}},$$
where
$$M(x, y)=\max\{d(x, y), d(x, Tx), d(y, Ty)\}.$$
\end{definition}
\begin{theorem}\label{Th3.4}
Let $(X, d, \preceq)$ be a complete partially ordered metric space with $x
\preceq y$ for all $x, y \in X$. Let $T: X\rightarrow X$ be a continuous mapping satisfying
\begin{description}
  \item [(1)]  $T$ is $\mathcal{Z}_{\mathcal{G}}$-Geraghty contraction;
  \item [(2)] $T$ is increasing;
  \item [(3)] there exists $x_{1}\in X$ such that $x_{1}\preceq Tx_{1}$.
\end{description}
Then $T$ has a fixed point $x^{*} \in X $ and $T$ is a Picard operator that
is, ${T^{n}x_{1}}$ converges to $x ^{*}$.

\begin{proof}
Define $\alpha:X\times X\rightarrow [0, +\infty)$ by
$$\alpha(x, y)=\begin{cases}
1, & ~ if ~ x\preceq y, \\
0, & otherwise.
\end{cases}$$
Since $T$ is a $Z_{\mathcal{G}}$-Geraghty contraction, we have
\begin{eqnarray*}
\mathcal{C}_{\mathcal{G}}&\leq& \zeta ( d(Tx , Ty), \beta (M(x, y)) M(x, y))\\
&<& \mathcal{G}(\beta (M(x, y)) M(x, y), d(Tx, Ty)).
\end{eqnarray*}
By definition of $\mathcal{G}$, we get
\begin{eqnarray*}
d(Tx, Ty)< \beta (M(x, y)) M(x, y),
\end{eqnarray*}
so,
\begin{eqnarray*}
\alpha(x, y) d(Tx, Ty)\leq d(Tx, Ty)<\beta(M(x, y)) M(x, y).
\end{eqnarray*}
Hence $T$ is $\mathcal{Z}_{(\alpha, \mathcal{G})}$-Geraghty contraction. Since $T$ is increasing, $\alpha(x, y)= 1$ implies  $\alpha(Tx, Ty)=
      1$ for all $x, y \in X$. Further if $z=Ty$, then $\alpha(y, z)=1$, this implies $\alpha(x, z)=1$. Thus, $T$ is triangular
$\alpha$-admissible. Condition $(2)$ implies that there exists $x_{1}\in X$ such that
       $\alpha(x_{1}, Tx_{1})=1$, and so condition $(3)$ of Theorem \ref{Th3.1} is
       satisfied.
Thus by Theorem \ref{Th3.1}, $T$ has a fixed point in $X$.
\end{proof}
\end{theorem}
Continuity of the mapping can be omitted in Theorem \ref{Th3.4} and fixed point result can be obtain with an extra condition given in the following theorem:
\begin{theorem}\label{Th3.5}
Let $(X, d, \preceq)$ be a complete partially ordered metric space with $x
\preceq y$ for all $x, y \in X$. Let $T:X\rightarrow X$ be a mapping satisfying
\begin{description}
  \item [(1)]  $T$ is $\mathcal{Z}_{\mathcal{G}}$-Geraghty contraction
      type mapping;
  \item [(2)] $T$ is increasing;
  \item [(3)] there exists $x_{1}\in X $ such that $x_{1}\preceq Tx_{1}$;
  \item [(4)] if $x_{n}$ is a sequence in $X$ such that $x_{n}\preceq Tx_{n+1}$ for all $n$ and $x_{n} \rightarrow x \in X$ as $n
      \rightarrow \infty $, then there exist a subsequence $\{x_{n(k)}\}$
      of $\{x_{n}\}$ such that $x_{n(k)}\preceq x$ for all $k$.
\end{description}
Then $T$ has a fixed point $x^{*} \in X $ and $T$ is a Picard operator that
is, ${T^{n}x_{1}}$ converges to $x ^{*}$.
\begin{proof}
Define $\alpha: X\times X\rightarrow [0, +\infty)$ by
$$\alpha(x, y)=\begin{cases}
1 & ~ if ~ x\preceq y, \\
0 & otherwise.
\end{cases}$$
Since $T$ is a $\mathcal{Z}_{\mathcal{G}}$-Geraghty contraction, we have
\begin{eqnarray*}
\mathcal{C}_{\mathcal{G}}&\leq& \zeta ( d(Tx, Ty), \beta (M(x, y)) M(x, y))\\
&<& \mathcal{G}(\beta (M(x, y)) M(x, y), d(Tx, Ty)).
\end{eqnarray*}
By the definition of $\mathcal{G}$, we obtain
\begin{eqnarray*}
d(Tx, Ty)< \beta (M(x, y)) M(x, y).
\end{eqnarray*}
This implies
\begin{eqnarray*}
\alpha(x, y) d(Tx, Ty)\leq d(Tx, Ty)< \beta (M(x, y)) M(x, y).
\end{eqnarray*}
Hence $T$ is generalized $\mathcal{Z}_{(\alpha, \mathcal{G})}$-Geraghty contraction. Since $T$ is increasing, $\alpha(x, y)= 1$ implies  $\alpha(Tx, Ty)=
1$ for all $x, y\in X$. Further if $z=Ty$, then $\alpha(y, z)=1$, this
implies $\alpha(x, z)=1$. Thus, $T$ is triangular $\alpha$-admissible.
      Condition $(2)$ implies that there exists $x_{1}\in X$ such that
       $\alpha(x, Tx)=1$ and so condition $(3)$ of Theorem \ref{Th3.2} is
       satisfied. Condition $(4)$ implies that the condition $(4)$ of Theorem \ref{Th3.2} is satisfied.
       Thus, all the conditions of Theorem \ref{Th3.2} are satisfied. Hence $T$ has a fixed point in $X$.
\end{proof}
\end{theorem}
\begin{remark}
Uniqueness of fixed point follows from Theorem \ref{Th3.4} (respectively Theorem \ref{Th3.5}) with the condition \\
{\bf $\mathcal{U}$:} For all $x,y\in Fix(T)$ with $x\preceq y$, there exists $z\in X$ such that $\alpha(x, z)\geq 1$ and
$\alpha(y, z)\geq 1$.
\end{remark}
\section{Application to Differential Equations}
Denote by $C([0, 1])$ the set of all continuous functions
defined on $[0, 1]$ and let $d :C([0, 1])\times C([0, 1])\to \mathbb{R}$ be
defined by
\begin{equation}\label{1}
d (x,y)=||x-y||_{\infty}=\max_{t\in[0,1]}|x(t)-y(t)|.
\end{equation}
It is well known that $(C([0, 1]),d )$ is a complete metric space.
Let us consider the two-point boundary value problem of the
second-order differential equation:
\begin{eqnarray}\label{2}
\notag
&&- \frac{d^2 x}{dt^2}=f(t,x(t)),~~t\in[0,1];\\
&&x(0)=x(1)=1,
\end{eqnarray}
where $f: [0, 1]\times \mathbb{R}\to \mathbb{R}$ is continuous.The Green function
associated to \eqref{2} is defined by
$$
G(t,s)=\left\{ \begin{array}{cc}
           t(1-s)&\text{if}~0\leq t\leq s\leq1,\\
           s(1-t)&\text{if}~0\leq s\leq t\leq1.
\end{array}
\right.
$$
Assume that the following conditions hold:
\begin{enumerate}
  \item[(i)] there exist a function $\xi:\mathbb{R}\times\mathbb{R}\to\mathbb{R}$  such that $$|f(t,a)-f(t,b)| \leq \max{\left\{|a-b|,|a-Ta|,|b-Tb|\right\}}$$
 for all $t\in[0,1]$ and $a,b\in \mathbb{R}$ with $\xi(a,b)>0$, where $T:C[0,1]\to C[0,1]$
is defined by $$Tx(t)=\int_{0}^{1}G(t,s)f(s,x(s))ds;$$
  \item[(ii)] there exists $x_0\in C[0,1]$ such that $\xi(x_0(t),Tx_0(t))\geq0$ for all $t\in[0,1]$;
  \item[(iii)] for each $t\in [0,1]$ and $x,y\in C[0,1]$, $\xi(x(t),y(t))>0$ implies $\xi(Tx(t),Ty(t))>0$;
  \item[(iv)] for each $t\in [0,1]$, if $\{x_n\}$ is e sequence in $C[0,1]$ such that $x_n\to x$ in $C[0,1]$ and $\xi(x_{n}(t),x_{n+1}(t))>0$ for all $n\in \mathbb{N}$, then $\xi(x_n(t),x(t))>0$ for all $n\in\mathbb{N}$.
\end{enumerate}
We now prove that existence of a solution of the mentioned
second-order differential equation.
\begin{theorem}
 Under assumptions $(i)-(iv)$, \eqref{2} has a solution
in $C^2([0, 1])$.
\begin{proof}
 It is well known that $x\in C^2([0, 1])$ is a solution of
\eqref{2} is equivalent to $x\in C([0, 1])$ is a solution of the integral
equation (see \cite{samet})
\begin{equation}\label{3}
x(t)=\int_{0}^{1}G(t,s)f(s,x(s))ds, ~t\in [0,1].
\end{equation}
Let $T:C[0,1]\to C[0,1]$ be a mapping defined by
\begin{equation}\label{4}
Tx(t)=\int_{0}^{1}G(t,s)f(s,x(s))ds.
\end{equation}
Suppose that $x,y\in C([0, 1])$ such that $\xi(x(t),y(t))\geq 0$ for all $t\in[0,1]$. By applying (i), we obtain that
\begin{eqnarray*}
&&|Tu(x)-Tv(x)|\\
&=&\int_{0}^{1}G(t,s)f(s,x(s))ds-\int_{0}^{1}G(t,s)f(s,y(s))ds\\
&=&\int_{0}^{1}G(t,s)[f(s,x(s))-f(s,y(s))]ds\\
&\leq&\Bigg(\int_{0}^{1}G(t,s)ds\Bigg)\Bigg(\int_{0}^{1}|f(s,x(s))-f(s,y(s))|ds\Bigg)\\
&\leq&\Bigg(\int_{0}^{1}G(t,s)ds\Bigg)\Bigg(\int_{0}^{1}(\max\{|x(s)-y(s)|,|x(s)-Tx(s)|,|y(s)-Ty(s)|\}ds\Bigg)\\
&\leq&\sup_{t\in[0,1]}\Bigg(\int_{0}^{1}G(t,s)ds\Bigg)\Bigg(\int_{0}^{1}(\max\left\{\sup_{s\in[0,1]}|x(s)-y(s)|,\right.\\
&&\left.\sup_{s\in[0,1]}|x(s)-Tx(s)|,\sup_{s\in[0,1]}|y(s)-Ty(s)|\right\}ds\Bigg)\\
&\leq&\sup_{t\in[0,1]}\Bigg(\int_{0}^{1}G(t,s)ds\Bigg)(\max\left\{||x-y||_{\infty},||x-Tx||_{\infty},||y-Ty||_{\infty}\right\})\int_{0}^{1}ds\\
&\leq&\sup_{t\in[0,1]}\Bigg(\int_{0}^{1}G(t,s)ds\Bigg)(M(x,y)).
\end{eqnarray*}
Since $\int_{0}^{1}G(t,s)ds = -(t^2/2) + (t/2)$, for all $t\in [0,1]$, we have\\
$\sup_{t\in[0,1]}(\int_{0}^{1}G(t,s)ds)=1/8$. It follows that
\begin{equation}\label{5}
||Tx-Ty||_{\infty}\leq\frac{1}{8}M(x,y).
\end{equation}
Let $\zeta(t, s)= \frac{1}{4}s -t,
\mathcal{G}(s, t)=s-t$ for all $s, t\in[0, \infty)$, $\mathcal{C}(\mathcal{G})=0$ and
$\beta(t)=\frac{1}{2}$ for all $t\geq 0$. Then it is clear that $\beta \in
\mathcal{F}$. Also define
$$
\alpha(x,y)=\left\{ \begin{array}{cc}
           1 &\text{if}~\xi(x(t),y(t))>0,t\in[0,1],\\
           0 &\text{otherwise}.
\end{array}
\right.
$$
Now
\begin{eqnarray}\label{6}
\notag
 \zeta(\alpha(x, y)d(Tx, Ty), \beta(M(x, y))M(x, y))&=&\frac{1}{4}\beta(M(x, y))M(x, y)- \alpha(x, y)d(Tx, Ty)\\
&=&\frac{1}{8}M(x, y)-d(Tx, Ty),
\end{eqnarray}
Then from \eqref{5}
$$\zeta(\alpha(x, y)d(Tx, Ty), \beta(M(x, y))M(x, y))\geq 0.$$
Therefore the mapping $T$ is a $\mathcal{Z}_{(\alpha,
\mathcal{G})}$-Geraghty contraction.

From (ii) there exists $x_0\in C[0,1]$ such that $\alpha(x_0,Tx_0)\geq 1$. Next by using (iii), we get the following assertions holding for all $x,y\in C[0,1]$
\begin{eqnarray*}
\alpha(x,y)\geq 1 &\Rightarrow& \xi(x(t),y(t))>0~~\text{for all}~t\in [0,1]\\
                  &\Rightarrow& \xi(Tx(t),Ty(t))>0~~\text{for all}~t\in [0,1]\\
                  &\Rightarrow& \alpha(Tx,Ty)\geq 1,
\end{eqnarray*}
hence $T$ is $\alpha$-admissible.

Applying Theorem \eqref{Th3.1}, we obtain that $T$ has a fixed point in
$C([0, 1])$; say $x$. Hence, $x$ is a solution of \eqref{2}.
\end{proof}
\end{theorem}

\subsection*{Acknowledgments}
This paper was funded by the University of Sargodha, Sargodha funded research project No. UOS/ORIC/2016/54. The first author,
therefore, acknowledges with thanks UOS for financial support.

\end{document}